\documentclass[11pt,a4paper]{article}
\usepackage{amsmath,amssymb,amsthm}
\usepackage{amscd}

\newtheorem{theorem}{Theorem}[section]
\newtheorem{lemma}[theorem]{Lemma}
\newtheorem{corollary}[theorem]{Corollary}
\newtheorem{proposition}[theorem]{Proposition}

{\theoremstyle{definition}
\newtheorem{condition}[theorem]{Condition}

\newtheorem{remark}[theorem]{Remark}

\newcommand{\bim}[3]{\mathord{\raisebox{-0.4ex}[0ex][0ex]{\scriptsize $#1$}{#2}\hspace{-0.2ex}\raisebox{-0.4ex}[0ex][0ex]{\scriptsize $#3$}}}

\newcommand{\lmo}[2]{\mathord{\raisebox{-0.4ex}[0ex][0ex]{\scriptsize $#1$}{#2}}}

\newcommand{\cG}{\mathcal{G}}
\newcommand{\recht}{\rightarrow}
\newcommand{\rL}{\operatorname{L}}
\newcommand{\actson}{\curvearrowright}
\newcommand{\cI}{\mathcal{I}}
\newcommand{\cC}{\mathcal{C}}
\newcommand{\rZ}{\operatorname{Z}}
\newcommand{\M}{\operatorname{M}}
\newcommand{\C}{\mathbb{C}}
\newcommand{\ot}{\otimes}
\newcommand{\Tr}{\operatorname{Tr}}
\newcommand{\cU}{\mathcal{U}}
\newcommand{\Aut}{\operatorname{Aut}}
\newcommand{\cH}{\mathcal{H}}
\newcommand{\cK}{\mathcal{K}}
\newcommand{\cL}{\mathcal{L}}
\newcommand{\si}{\sigma}
\newcommand{\diml}{\operatorname{dim}_{\ell}}
\newcommand{\dimr}{\operatorname{dim}_{r}}
\newcommand{\Comm}{\operatorname{Comm}}
\newcommand{\Bim}{\operatorname{Bim}}

\newcommand{\Bimod}{\operatorname{Bimod}}
\newcommand{\B}{\operatorname{B}}
\newcommand{\Ad}{\operatorname{Ad}}
\newcommand{\al}{\alpha}

\newcommand{\N}{\mathbb{N}}
\newcommand{\Hecke}{\operatorname{Hecke}}
\newcommand{\om}{\omega}
\newcommand{\D}{\operatorname{D}}
\newcommand{\pitil}{\widetilde{\pi}}
\newcommand{\Ker}{\operatorname{Ker}}
\newcommand{\vphi}{\varphi}
\newcommand{\cR}{\mathcal{R}}
\newcommand{\cZ}{\mathcal{Z}}
\newcommand{\Stab}{\operatorname{Stab}}
\newcommand{\PU}{\operatorname{P\mathcal{U}}}
\newcommand{\omtil}{\widetilde{\omega}}
\newcommand{\SL}{\operatorname{SL}}
\newcommand{\R}{\mathbb{R}}
\newcommand{\lspan}{\operatorname{span}}
\newcommand{\Q}{\mathbb{Q}}
\newcommand{\Z}{\mathbb{Z}}
\newcommand{\Gal}{\operatorname{Gal}}
\newcommand{\id}{\mathord{\operatorname{id}}}
\newcommand{\GL}{\operatorname{GL}}
\newcommand{\cP}{\mathcal{P}}
\newcommand{\dpr}{^{\prime\prime}}
\newcommand{\Om}{\Omega}
\newcommand{\F}{\mathbb{F}}
\newcommand{\cV}{\mathcal{V}}
\newcommand{\PSL}{\operatorname{PSL}}
\newcommand{\cN}{\mathcal{N}}
\newcommand{\rP}{\operatorname{P}}
\newcommand{\T}{\mathbb{T}}

\begin{document}

\begin{center}
{\LARGE\bf A classification of all finite index subfactors
 \vspace{0.5ex}\\ for a class of group-measure space II$_1$ factors}

\bigskip

{\sc by Steven Deprez$^{(1,3)}$\setcounter{footnote}{1}\footnotetext{Research assistant of the Research Foundation --
    Flanders (FWO)} and Stefaan Vaes$^{(2,3)}$\setcounter{footnote}{2}\footnotetext{Partially
    supported by ERC Starting Grant VNALG-200749, Research
    Programme G.0231.07 of the Research Foundation --
    Flanders (FWO) and K.U.Leuven BOF research grant OT/08/032.}\setcounter{footnote}{3}\footnotetext{Department of Mathematics;
    K.U.Leuven; Celestijnenlaan 200B; B--3001 Leuven (Belgium).
    \\ E-mails: steven.deprez@wis.kuleuven.be and stefaan.vaes@wis.kuleuven.be}}
\end{center}

\bigskip

\begin{abstract}
\noindent We provide a family of group measure space II$_1$ factors for which all finite index subfactors can be explicitly listed. In particular, the set of all indices of irreducible subfactors can be computed. Concrete examples show that this index set can be any set of natural numbers that is closed under taking divisors.
\end{abstract}

\section{Introduction and statement of main results}

Recall that the \emph{Jones index} \cite{Jo83} of an inclusion of II$_1$ factors $N \subset M$ is defined as the Murray-von Neumann dimension $[M:N]:=\dim(\rL^2(M)_N)$. The astonishing main result of \cite{Jo83} says that the index $[M:N]$ can only take values in $\cI \cup \{+\infty\}$ where $\cI:= \{4 \cos(\pi/n)^2 \mid n=3,4,5,\ldots\} \cup [4,+\infty)$. Conversely, Jones showed that all these values do arise as the index of a subfactor of the hyperfinite II$_1$ factor $R$. For general II$_1$ factors, Jones defines
$$\cI(M) := \{[M:N] \mid N \subset M \;\;\text{a finite index subfactor}\;\} \; .$$
For concrete II$_1$ factors $M$, determining $\cI(M)$ is extremely hard. Popa's deformation/rigidity theory (see \cite{Po06b} for a survey) has made it possible to compute invariants for several families of II$_1$ factors. This has been successfully applied for the fundamental group (see e.g.\ \cite{Po01,Po03,PV08a}) and the outer automorphism group (see e.g.\ \cite{IPP05,PV06,Va07,FV07}). In this paper, we apply Popa's theory to provide computations of $\cI(M)$.

Although it is known that $\cI(R) = \cI$, it is a major open problem to compute the irreducible counterpart $\cC(R)$, defined for arbitrary II$_1$ factors $M$ as
$$\cC(M) := \{[M:N] \mid N \subset M \;\;\text{an irreducible finite index subfactor}\;\} \;.$$
Recall that a subfactor $N \subset M$ is said to be \emph{irreducible} if $N' \cap M = \C1$.

In \cite{Va06,Va07}, examples of II$_1$ factors $M$ without non-trivial finite index subfactors were given. For these examples, $\cC(M) = \{1\}$ and $\cI(M) = \{n^2 \mid n \in \N\}$.

In this paper, we produce concrete II$_1$ factors $M$ for which $\cC(M)$ and $\cI(M)$ can be computed. In particular, given any set $\cP_0$ of prime numbers, we provide examples where $\cI(M)$ consists of the positive integers having all prime divisors in $\cP_0$. We prove that $\cC(M)$ ranges over all sets of natural numbers with the property of being closed under taking divisors.
Our examples $M$ are group measure space II$_1$ factors $\rL^\infty(X,\mu) \rtimes \Gamma$ associated with a free, ergodic, probability measure preserving (p.m.p.) action $\Gamma \actson (X,\mu)$ satisfying a number of conditions. To explain these conditions, we need the following concepts.

\emph{Cocycle superrigidity.} Given a p.m.p.\ action $\Gamma \actson (X,\mu)$ and a Polish group $\cG$, a $1$-cocycle for $\Gamma \actson (X,\mu)$ with values in $\cG$ is a Borel map $\om : \Gamma \times X \recht \cG$ satisfying $\om(gh,x) = \om(g, h \cdot x) \om(h,x)$ almost everywhere. A $1$-cocycle is said to be cohomologous to a group morphism, if there exists a group morphism $\delta : \Gamma \recht \cG$ and a measurable map $\vphi : X \recht \cG$ such that $\om(g,x) = \vphi(g \cdot x) \delta(g) \vphi(x)^{-1}$ almost everywhere. A p.m.p.\ action $\Gamma \actson (X,\mu)$ is said to be \emph{$\cG$-cocycle superrigid} if every $1$-cocycle with values in $\cG$ is cohomologous to a group morphism.

\emph{Bimodules.} A $P$-$Q$-bimodule (or correspondence) between von Neumann algebras $P$ and $Q$ is a Hilbert space $\cH$ equipped with a normal representation of $P$ and a normal anti-representation of $Q$ having commuting images. If $P$ and $Q$ are tracial, a bimodule $\bim{P}{\cH}{Q}$ is said to be of \emph{finite index} if the Murray-von Neumann dimensions $\dim(\cH_Q)$ and $\dim(\lmo{P}{\cH})$ are both finite.

\emph{Weak mixing.} A p.m.p.\ action $\Gamma \actson (X,\mu)$ is called \emph{weakly mixing} if $\C 1$ is the only non-zero finite dimensional $\Gamma$-invariant subspace of $\rL^2(X)$. Equivalently, $\Gamma \actson (X,\mu)$ is weakly mixing iff there exists a sequence $g_n$ in $\Gamma$ such that $\mu(g_n \cdot \cU \cap \cV) \recht \mu(\cU) \mu(\cV)$ for all measurable subsets $\cU,\cV \subset X$.

\begin{condition}\label{cond}
We say that the action $\Gamma \actson (X,\mu)$ satisfies Condition \ref{cond} if $\Gamma \actson (X,\mu)$ is free, weakly mixing and p.m.p.\ and if, denoting $A = \rL^\infty(X)$ and $M = A \rtimes \Gamma$, the following conditions are satisfied.
\begin{itemize}
\item Every non-zero finite index $M$-$M$-bimodule admits a non-zero finite index $A$-$A$-subbimodule.
\item If $\cG$ is a countable or a compact second countable group and if $\Gamma_1 < \Gamma$ is a finite index subgroup, then every $1$-cocycle for $\Gamma_1 \actson (X,\mu)$ with values in $\cG$ is cohomologous to a group morphism $\Gamma_1 \recht \cG$.
\end{itemize}
\end{condition}

When $\Gamma \actson (X,\mu)$ satisfies Condition \ref{cond} and $M = \rL^\infty(X) \rtimes \Gamma$, we explicitly determine, up to unitary conjugacy, all finite index subfactors of $M$. In particular, we compute the invariants $\cI(M)$ and $\cC(M)$. When $N \subset M$ is a finite index subfactor, $\bim{N}{\rL^2(M)}{M}$ is a finite index $N$-$M$-bimodule. In general, we call II$_1$ factors $P$ and $Q$ \emph{commensurable} if there exists a non-zero finite index $P$-$Q$-bimodule. We determine, up to isomorphism, all II$_1$ factors that are commensurable with $M = \rL^\infty(X) \rtimes \Gamma$.

Thanks to Sorin Popa's deformation/rigidity theory, it is possible to give examples of group actions satisfying the strong conditions in \ref{cond}. Popa proved in \cite{Po05,Po06a} cocycle superrigidity theorems with arbitrary countable or arbitrary compact second countable target groups, providing many group actions that satisfy the second part of Condition \ref{cond}. On the other hand, the first part of Condition \ref{cond} can be considered as a strengthening of the assumption that every automorphism of $M$ preserves globally (and up to unitary conjugacy) the Cartan subalgebra $A$. Examples of group measure space II$_1$ factors satisfying such conditions were obtained in \cite{Po01,Po04,IPP05,PV06,Va07,OP07,PV09}.

More concretely, we get the following examples of group actions satisfying both conditions in \ref{cond}. All of them are (generalized) Bernoulli actions: given the action $\Gamma \actson I$ of $\Gamma$ on the countable set $I$ and given a probability space $(X_0,\mu_0)$, we consider $\Gamma \actson (X_0,\mu_0)^I$ given by $(g \cdot x)_i = x_{g^{-1} \cdot i}$. In all of the following examples, cocycle superrigidity is provided by Popa's \cite[Theorem 0.1 and Proposition 3.6]{Po05}.

\begin{itemize}
\item By \cite[Theorem 7.1]{PV09}, the generalized Bernoulli actions $\Gamma \actson (X_0,\mu_0)^I$ associated with the following $\Gamma \actson I$ satisfy Condition \ref{cond}~: $\Gamma = \SL(n,\Z) *_\Sigma (\Sigma \times \Lambda)$, where $n \geq 3$, $\Lambda$ is an arbitrary non-trivial group and $\Sigma \cong \Z$ is generated by a hyperbolic element $A \in \SL(n,\Z)$ such that $\Sigma \cdot i$ is infinite for all $i \in I$.

More generally, let $\Gamma = \Gamma_1 *_\Sigma \Gamma_2$ be an amalgamated free product, where $\Sigma$ is an infinite amenable group. Let $\Gamma \actson I$. Among other examples, \cite[Theorem 7.1]{PV09} implies that, under the following assumptions, $\Gamma \actson (X_0,\mu_0)^I$ satisfies Condition \ref{cond}.
 \begin{itemize}
 \item $\Gamma_1$ has a normal, non-amenable subgroup $H$ with the relative property (T) such that $H \cdot i$ is infinite for all $i \in I$.
 \item $\Sigma$ is a proper normal subgroup of $\Gamma_2$ and $\Sigma \cdot i$ is infinite for all $i \in I$.
 \item $\Gamma$ admits a subgroup $G$ of infinite index such that $g \Sigma g^{-1} \cap \Sigma$ is finite for all $g \in \Gamma - G$.
 \end{itemize}

\item By \cite[Theorem 2.2]{Va07}, the generalized Bernoulli actions $\Gamma \actson (X_0,\mu_0)^I$ associated with the following $\Gamma \actson I$ satisfy Condition \ref{cond}~: $\PSL(n,\Z) \actson \rP(\Q^n)$ for $n \geq 3$ and $\SL(n,\Z) \ltimes \Z^n \actson \Z^n$ for $n \geq 2$. For more examples, see \cite[Examples 2.5]{Va07}.
\end{itemize}

In order to state the main result of this paper, we need the following concepts.

{\it Stable isomorphism of II$_1$ factors.} Recall that $M^t$ denotes the amplification of a II$_1$ factor $M$ by $t > 0$~: up to isomorphism, $M^t = p(\M_n(\C) \ot M)p$ for some orthogonal projection $p \in \M_n(\C) \ot M$ satisfying $(\Tr \ot \tau)(p) = t$. We say that the II$_1$ factors $M$ and $N$ are \emph{stably isomorphic} if there exists $t>0$ and an isomorphism $\pi : N \recht M^t$. Associated with $\pi$, is the natural $M$-$N$-bimodule $\bim{M}{\cH^\pi}{N}$ given by $\bim{M}{(\M_{1,n}(\C) \ot \rL^2(M))p}{\,\pi(N)}$. We can equivalently define a \emph{stable isomorphism} between $M$ and $N$ as an $M$-$N$-bimodule $\bim{M}{\cK}{N}$ with the property that the right $N$-action on $\cK$ equals the commutant of the left $M$-action on $\cK$. Every stable isomorphism $\bim{M}{\cK}{N}$ is unitarily equivalent with $\bim{M}{\cH^\pi}{N}$ for an isomorphism $\pi : N \recht M^t$ that is uniquely determined up to unitary conjugacy.

{\it Commensurate subgroups and commensurators.} Recall that subgroups $\Gamma, G < \cG$ are called commensurate if $\Gamma \cap G$ has finite index in both $\Gamma$ and $G$. If $\Gamma < \cG$, the commensurator of $\Gamma$ inside $\cG$ is defined as the group of all $g \in \cG$ such that $g \Gamma g^{-1} \cap \Gamma$ has finite index in both $\Gamma$ and $g \Gamma g^{-1}$.

{\it Projective representations.} Recall that a map $\pi : G \recht \cU(\cK)$ from a countable group $G$ to the unitary group of a Hilbert space $\cK$ is called a projective representation if $\pi(g) \pi(h) = \Omega(g,h) \pi(gh)$ for some map $\Omega : G \times G \recht \T$. We call $\Omega$ the \emph{obstruction $2$-cocycle} of $\pi$. The group of $2$-cocycles $\Omega : G \times G \recht \T$ is denoted by $\rZ^2(G,\T)$.

{\it Cocycle crossed products.} Whenever $G \actson (X,\mu)$ and $\Omega \in \rZ^2(G,\T)$, one constructs the cocycle crossed product $\rL^\infty(X) \rtimes_\Omega G$, generated by $\rL^\infty(X)$ and unitaries $(u_g)_{g \in G}$ satisfying
$$u_g^* a u_g = a(g \, \cdot \,) \;\; , \;\; u_g u_h = \Omega(g,h)u_{gh} \;\;\text{and}\;\; \tau(a u_g) = \begin{cases} \int_X a \; d\mu &\;\;\text{if}\; g = e \; , \\ 0 &\;\;\text{if}\; g\neq e \; ,\end{cases}$$
for all $a \in \rL^\infty(X)$, and $g,h \in G$.

{\it Connes tensor product.} When $\bim{P}{\cH}{Q}$ and $\bim{Q}{\cK}{N}$ are bimodules between the tracial von Neumann algebras $P$, $Q$ and $N$, we denote by $\cH \ot_Q \cK$ the Connes tensor product, which is a $P$-$N$-bimodule.

\begin{theorem}\label{thm.main}
Suppose that the action $\Gamma \actson (X,\mu)$ satisfies Condition \ref{cond}. Put $A = \rL^\infty(X)$ and $M = A \rtimes \Gamma$.
\begin{enumerate}
\item Up to stable isomorphism, the II$_1$ factors that are commensurable with $M$ are precisely given as $A \rtimes_\Omega G$, where $G \subset \Aut(X,\mu)$ is commensurate with $\Gamma$ and $\Omega \in \rZ^2(G,\T)$ is a $2$-cocycle that arises from a finite dimensional projective representation of $G$.
\item With $(G,\Omega)$ and $(H,\omega)$ as in 1, put $P = A \rtimes_\Omega G$ and $Q = A \rtimes_\omega H$. Every irreducible finite index $P$-$Q$-bimodule $\bim{P}{\cH}{Q}$ is unitarily equivalent with
    $$\bim{P}{\cH(\gamma,\pi)}{Q} := \rL^2(P) \otimes_{P_0} \cK(\gamma,\pi) \otimes_{Q_0} \rL^2(Q) \; ,$$
    where $P_0 \subset P$ and $Q_0 \subset Q$ are the finite index subfactors defined by
    $$P_0 = A \rtimes_\Omega (G \cap \gamma H \gamma^{-1}) \quad\text{and}\quad Q_0 = A \rtimes_\omega (H \cap \gamma^{-1} G \gamma)$$
    for some $\gamma \in \Aut(X,\mu)$ in the commensurator of $\Gamma$, and where the bimodule $\bim{P_0}{\cK(\gamma,\pi)}{Q_0} = \bim{\psi_\pi(P_0)}{(\rL^2(Q_0) \ot K)}{Q_0}$ is given by an irreducible projective representation $\pi : G \cap \gamma H \gamma^{-1} \recht K$ with obstruction $2$-cocycle $\Omega \, \overline{\omega \circ \Ad \gamma^{-1}}$ and corresponding inclusion
    $$\psi_\pi : P_0 \recht Q_0 \ot \B(K) : \psi_\pi(a u_g)= a(\gamma \, \cdot \,) u_{\gamma^{-1} g \gamma} \otimes \pi(g) \; .$$
\end{enumerate}
\end{theorem}

\begin{remark}
In order to deduce from Theorem \ref{thm.main} an entirely explicit list of all commensurable II$_1$ factors, we need to know the commensurator of $\Gamma$ inside $\Aut(X,\mu)$. Usually, this is a hard problem. But, when $\Gamma \actson I$ has the property that $(\Stab i) \cdot j$ is infinite for all $i \neq j$, then the commensurator of $\Gamma$ inside $\Aut\bigl((X_0,\mu_0)^I\bigr)$ has been computed in \cite[Lemma 6.15]{PV09} (see also \cite[Proof of Theorem 5.4]{PV06}). It is generated by the following two types of elementary commensurating automorphisms. Firstly, whenever $\eta$ is a permutation of $I$ that commensurates $\Gamma$, we have $\Delta$ given by $\bigl(\Delta(x)\bigr)_i = x_{\eta^{-1}(i)}$. Secondly, for every orbit $\Gamma \cdot i$ and every $\Delta_0 \in \Aut(X_0,\mu_0)$, we have $\Delta$ given by $\bigl(\Delta(x)\bigr)_j = \Delta_0(x_j)$ for $j \in \Gamma \cdot i$ and $\bigl(\Delta(x)\bigr)_j = x_j$ for $j \not\in \Gamma \cdot i$.
\end{remark}

For a given II$_1$ factor $M$, the finite index $M$-$M$-bimodules form a C$^*$-tensor category $\Bimod(M)$. More generally, one can build the so-called \emph{C$^*$-bicategory} $\Comm(M)$ of commensurable II$_1$ factors~: the objects (or $0$-cells) in this category are the II$_1$ factors $N$ that are commensurable with $M$, the morphisms (or $1$-cells) from $N_1$ to $N_2$ are the finite index $N_1$-$N_2$-bimodules with composition given by the Connes tensor product and finally, the morphisms between two morphisms (or $2$-cells) are given by the bimodular bounded operators between two bimodules. In Section \ref{sec.categorical}, we reinterpret Theorem \ref{thm.main} and prove that for the II$_1$ factors $M$ given by Theorem \ref{thm.main}, the bicategory $\Comm(M)$ is equivalent with a bicategory $\Hecke(\Gamma < \cG)$ associated with the Hecke pair $\Gamma < \cG$, where $\cG$ denotes the commensurator of $\Gamma$ inside $\Aut(X,\mu)$.

In particular, we get an explicit description of $\Bimod(M)$ as C$^*$-tensor category; a problem that was left open in \cite{Va07}, although all finite index bimodules could be described up to unitary equivalence.

Whenever $N \subset M$ is a finite index subfactor, $\bim{M}{\rL^2(M)}{N}$ is a finite index bimodule of left dimension $1$ and right dimension $[M:N]$. Hence, whenever $\Gamma \actson (X,\mu)$ satisfies Condition \ref{cond}, Theorem \ref{thm.main} provides a complete description of all finite index subfactors of $\rL^\infty(X) \rtimes \Gamma$, up to unitary conjugacy. We can make this description more concrete in the following way.

\begin{corollary}
Suppose that the action $\Gamma \actson (X,\mu)$ satisfies Condition \ref{cond}. Put $A = \rL^\infty(X)$ and $M = A \rtimes \Gamma$.
Up to unitary conjugacy, all finite index subfactors of $M$ are provided by the following construction.

The data for the construction consists of a subgroup $G < \Aut(X,\mu)$ that is commensurate with $\Gamma$, a $2$-cocycle $\Omega \in \rZ^2(G,\T)$, elements $\gamma_1,\ldots,\gamma_n$ in the commensurator of $\Gamma$ inside $\Aut(X,\mu)$ and finite dimensional projective representations $\pi_i : \Gamma \cap \gamma_i G \gamma_i^{-1} \recht \cU(K_i)$ with obstruction $2$-cocycle $\overline{\Omega \circ \Ad \gamma_i^{-1}}$.

Given these data, put $\ell_i = [G:G \cap \gamma_i^{-1} \Gamma \gamma_i]$ and $r_i = [\Gamma : \Gamma \cap \gamma_i G \gamma_i^{-1}]$.
Amplifying
$$\begin{CD}
\bigl(A \rtimes_\Omega G\bigr)^{1/\ell_i} @>{\text{\rm tunnel construction}}>{\text{\rm inclusion with index}\;\;\ell_i}> A \rtimes_{\Om \circ \Ad \gamma_i^{-1}} (\Gamma \cap \gamma_i G \gamma_i^{-1}) @>{a(\,\cdot\,)u_g \mapsto a(\gamma_i \, \cdot) u_{\gamma_i^{-1} g \gamma_i}}>{\text{\rm inclusion with index}\;\;\ell_i}> A \rtimes_\Omega G \\
@. @V{\text{\rm inclusion with index}\;\; r_i \dim(\pi_i)^2}V{au_g \mapsto au_g \ot \pi(g)}V @. \\
@. (A \rtimes \Gamma) \ot \B(K_i) @.
\end{CD}$$
yields an inclusion $A \rtimes_\Omega G \recht M^{\ell_i \dim(\pi_i)}$ of index $\ell_i r_i \dim(\pi_i)^2$. Reducing the diagonal product of these inclusions, provides a subfactor
$$(A \rtimes_\Omega G)^{1/\ell} \subset M$$
of index $\ell r$, where $\ell = \sum_{i=1}^n \ell_i \dim(\pi_i)$ and $r = \sum_{i=1}^n r_i \dim(\pi_i)$.
\end{corollary}

When the group $\Gamma$ has no non-trivial finite dimensional unitary representations (and in particular, no non-trivial finite index subgroups), the formulation becomes easier and we get the following result. Concrete examples of the invariants $\cC(M)$ and $\cI(M)$ are provided by Corollary \ref{cor.concrete}, yielding II$_1$ factors $M$ such that $\cC(M)$ is any prescribed set of natural numbers with the property of being stable under taking divisors.

\begin{corollary}\label{cor.case-no-fin-dim-rep}
Suppose that the action $\Gamma \actson (X,\mu)$ satisfies Condition \ref{cond} and that the group $\Gamma$ has no non-trivial finite dimensional unitary representations. Put $M = \rL^\infty(X) \rtimes \Gamma$ and denote by $\cG$ the normalizer of $\Gamma$ inside $\Aut(X,\mu)$. Then,
\begin{align*}
\cC(M) &= \bigl\{ |G| \;\big|\; G \;\;\text{is a finite subgroup of}\;\; \cG/\Gamma \bigr\} \; , \\
\cI(M) &= \bigl\{ n^2 |G| \;\big|\; G \;\;\text{is a finite subgroup of}\;\; \cG/\Gamma \;\;\text{and}\;\; n \in \N \setminus \{0\} \bigr\}   \; .
\end{align*}
\end{corollary}

\section{Proof of Theorem \ref{thm.main} and its corollaries}

In order to make our proofs more readable, we split it up into a several independent lemmas.
We start with a more precise version of the first part of \cite[Theorem 6.4]{Va07}.

\begin{lemma} \label{lem.projrep}
Let $G \actson (X,\mu)$ be a free, p.m.p.\ action such that, for every finite index subgroup $G_1 < G$, the action $G_1 \actson (X,\mu)$ is ergodic and cocycle superrigid with arbitrary compact second countable target groups. Let $\Omega \in \rZ^2(G,\T)$ and put $A = \rL^\infty(X)$. Let $Q$ be a II$_1$ factor and $\psi : A \rtimes_\Omega G \recht Q$ an irreducible finite index embedding, with corresponding bimodule $\bim{\psi(A \rtimes_\Omega G)}{\rL^2(Q)}{Q}$.

There exists a finite index subgroup $G_1 < G$, a finite dimensional projective representation $\pi : G_1 \recht \cU(K)$ with obstruction $2$-cocycle $\Omega_\pi$, a projection $p \in Q \cap \psi(A)'$ of trace $([G : G_1] \, \dim \pi)^{-1}$ and, writing $\Omega_1 := \Omega \overline{\Omega_\pi}$, a finite index inclusion $\psi_1 : A \rtimes_{\Omega_1} G_1 \recht p Q p$ such that $\psi_1(a) = \psi(a)p$ for all $a \in A$, $\psi_1(A)$ is maximal abelian in $pQp$ and
$$\bim{\psi(A \rtimes_\Omega G)}{\rL^2(Q)}{Q} \cong \ \raisebox{-1.3ex}{\scriptsize $A \rtimes_\Omega G$}\Bigl(\rL^2(A \rtimes_\Omega G) \underset{A \rtimes_\Omega G_1}{\ot} (p \rL^2(Q) \ot K)\Bigr)\hspace{-0.5ex}\raisebox{-1.3ex}{\scriptsize $Q$}$$
with left module action of $au_g \in A \rtimes_{\Omega} G_1$ on $p \rL^2(Q) \ot K$ given by $\psi_1(au_g) \ot \pi(g)$.

Moreover, $p$ can be chosen such that $\lspan \{\psi(u_g) p \psi(u_g)^* \mid g \in G\}$ is finite dimensional.
\end{lemma}
\begin{proof}
Put $P = A \rtimes_\Omega G$. Since $\psi(P) \subset Q$ has finite index, also $\psi(P) \cap \psi(A)' \subset Q \cap \psi(A)'$ has finite index (see e.g.\ \cite[Lemma A.3]{Va07}). Put $B = Q \cap \psi(A)'$. Since $\psi(A) \subset B$ has finite index, $B$ is of finite type I. Also, $\psi(A) \subset \cZ(B)$ has finite index. Moreover, $(\Ad \psi(u_g))_{g \in G}$ defines an ergodic action of $G$ on $\cZ(B)$. It follows that we can take a $*$-isomorphism $\theta : \cZ(B) \recht \rL^\infty(X \times \{1,\ldots,m\})$ such that $\theta(\psi(a)) = a \otimes 1$ for all $a \in A$. Hence, $\theta$ conjugates the action $(\Ad \psi(u_g))_{g \in G}$ on $\cZ(B)$ with the action $G \actson X \times \{1,\ldots,m\}$ given by
$$g \cdot (x,i) = (g \cdot x, \om(g,x) i)$$ where $\om : G \times X \recht S_m$ is a $1$-cocycle with values in the symmetric group $S_m$. By cocycle superrigidity, we may assume that $G \actson \{1,\ldots,m\}$ transitively and that $g \cdot (x,i) = (g \cdot x, g \cdot i)$. Define $G_1 = \Stab 1$, $p_0 = \theta^{-1}\bigl(\chi_{X \times \{1\}}\bigr)$ and $P_0 = A \rtimes_\Omega G_1$. Define $\psi_0 : P_0 \recht p_0 Q p_0$ as $\psi_0(d) = p_0 \psi(d)$ for all $d \in P_0$. By construction,
$$\bim{\psi(P)}{\rL^2(Q)}{Q} \cong \rL^2(P) \ot_{P_0} \bigl(\bim{\psi(P_0)}{p_0 \rL^2(Q)}{Q}\bigr) \; .$$
Put $Q_0 = p_0 Q p_0$ and $B_0 = Q_0 \cap \psi_0(A)'$. By construction, $\cZ(B_0) = \psi_0(A)$ and $(\Ad \psi_0(u_g))_{g \in G_1}$ defines an ergodic action of $G_1$. Since $B_0$ is finite and of type I, we can take a finite dimensional Hilbert space $K$ and a $*$-isomorphism $\gamma : B_0 \recht A \ot \B(K)$ such that $\gamma(\psi(a)) = a \ot 1$ for all $a \in A$ and such that $\gamma$ conjugates the action $(\Ad \psi(u_g))_{g \in G}$ on $B_0$ with the action $(\al_g)_{g \in G_1}$ on $A \ot \B(K) = \rL^\infty(X,\B(K))$ given by
$$(\al_{g^{-1}}(a))(x) = \beta(g,x)^{-1} (a(g \cdot x)) \quad\text{for}\;\; a \in \rL^\infty(X,\B(K)), g \in G_1, x \in X$$
where $\beta : G_1 \times X \recht \Aut(\B(K)) \cong \PU(K)$ is a $1$-cocycle. By cocycle superrigidity, we may assume that $\beta(g,x) = \Ad \pi(g)$ for some projective representation $\pi : G_1 \recht \cU(K)$.

Define unitaries $v_g \in Q_0$ as $v_g = \psi_0(u_g) \gamma^{-1}(1 \ot \pi(g)^*)$. Put $\Omega_1 = \Omega \overline{\Omega_\pi} \in \rZ^2(G_1,\T)$. Let $q_1 \in \B(K)$ be a minimal projection and put $p_1 = \gamma^{-1} (1 \ot q_1)$. For all $g \in G_1 \setminus \{e\}$, we have $E_{B_0}(\psi_0(u_g)) = 0$. Hence, also $E_{B_0}(v_g) = 0$. By construction, $v_g$ commutes with $p_1$ and we obtain a well defined finite index inclusion
$$\psi_1 : A \rtimes_{\Omega 1} G_1 \recht p_1 Q p_1 : \psi_1(a u_g) = a v_g p_1 \; .$$
This inclusion satisfies all desired properties.
\end{proof}

The next lemma is almost literally contained in \cite[Lemma 7.3 and proof of Theorem 7.1]{PV09}. We repeat it here for the convenience of the reader. We call a measurable map $\Delta : X \recht Y$ between probability spaces $(X,\mu)$ and $(Y,\eta)$, a \emph{local isomorphism,} if, up to measure zero, we can partition $X$ into non-negligible subsets $X_n$, $n \in \N$, such that the restriction $\Delta|_{X_n}$ is a non-singular isomorphism between $X_n$ and a non-negligible subset of $Y$. If moreover all $\Delta|_{X_n}$ are measure preserving, we call $\Delta$ a local m.p.\ isomorphism.

\begin{lemma}\label{lem.quotient}
Let $G \actson (X,\mu)$ be a free, p.m.p.\ action. Assume that all finite index subgroups of $G$ act ergodically on $(X,\mu)$ and that $G \actson (X,\mu)$ is cocycle superrigid with arbitrary countable target groups. Let $H \actson (Y,\eta)$ be an arbitrary free, ergodic, p.m.p.\ action.

If $\Delta : X \recht Y$ is a local isomorphism satisfying $\Delta(G \cdot x) \subset H \cdot \Delta(x)$ for almost all $x \in X$, there exists
\begin{itemize}
\item a finite group $\Lambda$ acting freely on $(X,\mu)$ and satisfying $g \Lambda g^{-1} = \Lambda$ for all $g \in G$,
\item an isomorphism of probability spaces $\Psi : X/\Lambda \recht Y$,
\item a group homomorphism $\delta : G \recht H$ with $\Ker \delta = G \cap \Lambda$,
\end{itemize}
such that $\psi(g \cdot x) = \delta(g) \cdot \psi(x)$ and $\Delta(x) \in H \cdot \Psi(x)$ for all $g \in G$ and almost all $x \in X$.
\end{lemma}
\begin{proof}
Define $\om : G \times X \recht H$ such that $\Delta(g \cdot x) = \om(g,x) \cdot \Delta(x)$ almost everywhere. Then, $\om$ is a $1$-cocycle and cocycle superrigidity provides a measurable map $\vphi : X \recht H$ such that, defining $\Psi(x) := \vphi(x) \cdot \Delta(x)$, we have $\Psi(g \cdot x) = \delta(g) \cdot \Psi(x)$ almost everywhere, for some group morphism $\delta : G \recht H$. By construction, $\Psi$ is still a local isomorphism.

Define the equivalence relation $\cR = \{(x,y) \in X \times X \mid \Psi(x) = \Psi(y)\}$. Since $\Psi$ is a local isomorphism and $\mu(X) < \infty$, almost every $x \in X$ has a finite $\cR$-equivalence class. Moreover, the function $x \mapsto \#\{y \in X \mid x\cR y\}$ is $G$-invariant and hence, almost everywhere equal to a constant that we denote by $m$. Denote by $Y_0 \subset Y$ the essential range of $\Psi$. It follows that, up to measure zero, we can partition $X$ into $X_1,\ldots,X_m$ such that $\Psi_i := \Psi|_{X_i}$ is an isomorphism between $X_i$ and $Y_0$. By construction, all $\Psi_i$ scale the involved measures by the same constant. Hence, $T_i := \Psi_i^{-1} \circ \Psi$ is a local m.p.\ isomorphism. Denote by $S_m$ the permutation group of $\{1,\ldots,m\}$. The formula
$$T_i(g \cdot x) = g \cdot T_{\eta(g,x)^{-1} i}(x)$$
defines a $1$-cocycle $\eta : G \times X \recht S_m$. By cocycle superrigidity, we find a measurable $\gamma : X \recht S_m$ and a group morphism $\pi : G \recht S_m$ such that, writing $R_i(x) = T_{\gamma(x)i}(x)$, we have $R_i(g \cdot x) = g \cdot R_{\pi(g)^{-1} i}(x)$. By construction, $R_i$ is a local m.p.\ isomorphism. Since the essential range of $R_i$ is globally invariant under the finite index subgroup $\Ker \pi < G$, we conclude that $R_i \in \Aut(X,\mu)$. By construction, the equivalence relation $\cR$ is the union of the graphs of $R_1,\ldots,R_m$. Since almost every $\cR$-equivalence class has $m$ elements, we deduce that the graphs of the $R_i$ are essentially disjoint. For all $i,j$, we find $k$ such that $\{x \in X \mid R_i(R_j(x)) = R_k(x)\}$ is non-negligible. But this last set is globally invariant under $\Ker \pi$, showing that $R_i \circ R_j = R_k$ almost everywhere. Hence, $\Lambda := \{R_1,\ldots,R_m\}$ is a subgroup of $\Aut(X,\mu)$, that defines an essentially free action $\Lambda \actson (X,\mu)$ and that satisfies all desired properties.
\end{proof}

\begin{lemma}\label{lem.noquotient}
Let $G \actson (X,\mu)$ be a free, weakly mixing, p.m.p.\ action. Let $\Lambda$ be a non-trivial finite group acting freely and p.m.p.\ on $(X,\mu)$. Assume that $g \Lambda g^{-1} = \Lambda$ for all $g \in G$. Then, $G$ admits a finite index subgroup $G_1$ such that the action
$$\frac{G_1}{G_1 \cap \Lambda} \actson \frac{X}{\Lambda}$$
is not cocycle superrigid with arbitrary countable target groups.
\end{lemma}
\begin{proof}
Take $G_1 < G$ of finite index such that $G_1$ and $\Lambda$ commute inside $\Aut(X,\mu)$. Define $\cG$ as the subgroup of $\Aut(X,\mu)$ generated by $G_1$ and $\Lambda$. Since $\Lambda$ is finite and $\Lambda \actson (X,\mu)$ is free, choose a measurable $\psi : X \recht \Lambda$ such that $\psi(\lambda \cdot x) = \lambda \psi(x)$ for all $\lambda \in \Lambda$ and almost all $x \in X$. Define
$$\om : G_1 \times X \recht \cG : \om(g,x) = \psi(g \cdot x)^{-1} g \psi(x) \; .$$
Since $G_1$ and $\Lambda$ commute, it is easy to check that $\om$ is actually a $1$-cocycle
$$\om : \frac{G_1}{G_1 \cap \Lambda} \times \frac{X}{\Lambda} \recht \cG \; .$$
If this $1$-cocycle were to be cohomologous to a group morphism, we would find a group morphism $\delta : G_1 \recht \cG$ and a $\Lambda$-invariant measurable map $\vphi : X \recht \cG$ such that
$$\psi(g \cdot x)^{-1} g \psi(x) = \vphi(g \cdot x) \delta(g) \vphi(x)^{-1} \; .$$
Define $F(x) = \psi(x) \vphi(x)$. It follows that $F(g \cdot x) = g F(x) \delta(g)^{-1}$. Hence, the set
$$\cU = \{(x,y) \in X \times X \mid F(x) = F(y)\}$$
is invariant under the diagonal $G_1$-action. Since $F$ takes only countably many values, $\cU$ is non-negligible. By weak mixing of $G \actson (X,\mu)$, it follows that $\cU$ has complement of measure zero. Hence, $F$ is constant almost everywhere. But then, $\psi$ follows $\Lambda$-invariant, which is a contradiction.
\end{proof}

\begin{lemma}\label{lem.untwist-pseudo}
Let $G \actson (X,\mu)$ be a free, weakly mixing, p.m.p.\ action that is cocycle superrigid with compact second countable target groups. Assume that $\Omega \in \rZ^2(G,\T)$ is the obstruction $2$-cocycle of a finite dimensional projective representation and that $\om : G \times X \recht \T$ is a measurable map satisfying
$$\om(gh,x) = \Omega(g,h) \, \om(g,h \cdot x) \, \om(h,x)$$
almost everywhere.

Then, there exists a measurable map $\vphi : X \recht \T$ and a map $\rho : G \recht \T$ such that
$$\om(g,x) = \overline{\vphi(g \cdot x)} \, \rho(g) \, \vphi(x)$$
almost everywhere. In particular, $\Omega$ is cohomologous to the trivial $2$-cocycle.
\end{lemma}
\begin{proof}
Assume that $\Omega = \Omega_\pi$ for a finite dimensional projective representation $\pi : G \recht \cU(K)$. Consider $\cU(K)$ as a compact group and define $\omtil : G \times X \recht \cU(K)$ by $\omtil(g,x) = \om(g,x) \pi(g)$. Then, $\om$ is a $1$-cocycle and cocycle superrigidity gives a measurable map $\vphi : X \recht \cU(K)$ and a homomorphism $\rho_1 : G \recht \cU(K)$ such that
$$\om(g,x) \pi(g) = \vphi(g \cdot x)^* \rho_1(g) \vphi(x) \; .$$
In the quotient $\PU(K)$, we get the equality $\vphi(g \cdot x) = \rho_1(g) \vphi(x) \pi(g)^*$ almost everywhere. Weak mixing of $G \actson (X,\mu)$ implies that $\vphi$ is almost everywhere constant in $\PU(K)$ (see e.g.\ \cite[Lemma 5.4]{PV08b}). Replacing $\vphi(\,\cdot\,)$ by $u \vphi(\,\cdot\,)$ for the appropriate $u \in \cU(K)$ and replacing $\rho_1$ by $(\Ad u^*) \circ \rho_1$, we may assume that $\vphi$ takes values in $\T \subset \cU(K)$. Then, $\rho_1(g) = \rho(g) \pi(g)$ for some map $\rho : G \recht \T$, proving the lemma.
\end{proof}

Fix a group action $\Gamma \actson (X,\mu)$ satisfying Condition \ref{cond}. Put $A = \rL^\infty(X)$ and $M = A \rtimes \Gamma$.

Define $\cG \subset \Aut(X,\mu)$ as the commensurator of $\Gamma$ inside $\Aut(X,\mu)$. By \cite[Lemma 6.11]{Va07}, every $\gamma \in \cG \setminus \{e\}$ acts essentially freely on $(X,\mu)$. So, whenever $G < \cG$ is a subgroup commensurate with $\Gamma$, the action $G \actson (X,\mu)$ is free and weakly mixing.

Whenever $G < \cG$ is commensurate with $\Gamma$, the action $G \actson (X,\mu)$ is cocycle superrigid with countable or with compact second countable groups. Indeed, by the assumption in Condition \ref{cond}, we have cocycle superrigidity of $G \cap \Gamma$ acting on $(X,\mu)$. Since $G \cap \Gamma$ acts weakly mixingly, cocycle superrigidity of $G \actson (X,\mu)$ follows from \cite[Proposition 3.6]{Po05}.

If $G < \cG$ is commensurate with $\Gamma$ and if $\Omega \in \rZ(G,\T)$ is the obstruction $2$-cocycle of the finite dimensional projective representation $\pi$ on $K$, the II$_1$ factor $A \rtimes_\Omega G$ is commensurable with $M$. Indeed, the embedding $A \rtimes_\Omega G \recht (A \rtimes G) \otimes \B(K) : au_g \mapsto au_g \ot \pi(g)$ proves the commensurability of $A \rtimes_\Omega G$ and $A \rtimes G$. The latter is commensurate with its finite index subfactor $A \rtimes (G \cap \Gamma)$, which is in turn commensurate with $A \rtimes \Gamma = M$.

\begin{proof}[{\bf Proof of Theorem \ref{thm.main}, point 2.}]
Let $G,H < \cG$ be subgroups that are commensurate with $\Gamma$. Take $\Omega \in \rZ^2(G,\T)$ and $\om \in \rZ^2(H,\T)$, both being obstruction $2$-cocycles for finite dimensional projective representations. Put $P = A \rtimes_\Omega G$ and $Q = A \rtimes_\omega H$. Let $\bim{P}{\cH}{Q}$ be an irreducible finite index $P$-$Q$-bimodule.

We first prove that $\cH$ is a direct sum of finite index $A$-$A$-subbimodules. In the paragraph preceding this proof, we constructed non-zero finite index bimodules $\bim{M}{\cK}{P}$ and $\bim{Q}{\cL}{M}$. By construction, $\cK$ and $\cL$ can be taken as a direct sum of finite index $A$-$A$-subbimodules. Consider now the $M$-$M$-bimodule $\cH' := \cK \ot_P \cH \ot_Q \cL$.

Condition \ref{cond} implies that every finite index $M$-$M$-bimodule $\cH'$ is a direct sum of finite index $A$-$A$-subbimodules. Indeed, it suffices to consider irreducible $\cH'$. By Condition \ref{cond}, we find a non-zero finite index $A$-$A$-subbimodule $\cH^{\prime\prime} \subset \cH'$. But then, all $u_g \cdot \cH^{\prime\prime} \cdot u_h$, $g,h \in \Gamma$ are finite index $A$-$A$-subbimodules. By irreducibility of $\cH'$, they densely span $\cH'$.

So, $\cH'$ is a direct sum of finite index $A$-$A$-subbimodules. Then, the same is true for the $P$-$Q$-bimodule $\overline{\cK} \ot_M \cH' \ot_M \overline{\cL}$. By construction, this last $P$-$Q$-bimodule contains $\bim{P}{\cH}{Q}$. We have proved that $\cH$ is a direct sum of finite index $A$-$A$-subbimodules.

Take the irreducible finite index inclusion $\psi_0 : P \recht p (\M_m(\C) \ot Q) p$ such that $\bim{P}{\cH}{Q}$ is isomorphic with $\bim{\psi_0(P)}{p (\C^m \ot \rL^2(Q))}{Q}$. We denote by $\D_m \subset \M_m(\C)$ the subalgebra of diagonal matrices. By \cite[Lemma 6.5]{Va07}, we may assume that $\psi_0(A) \subset (\D_m \ot A)p$.

Lemma \ref{lem.projrep} yields now the following data:
\begin{itemize}
\item a finite index subgroup $G_1 < G$,
\item a finite dimensional projective representation $\pi : G_1 \recht \cU(K)$ with obstruction $2$-cocycle $\Omega_\pi$,
\item using the notations $\Omega_1 = \Omega \overline{\Omega_\pi}$ and $P_1 = A \rtimes_{\Omega_1} G_1$, a finite index inclusion
$\psi : P_1 \recht q(\M_n(\C) \ot Q) q$ with $\psi(A) = (\D_n \ot A)q$ for some projection $q \in D_n \ot A$,
\end{itemize}
such that $\bim{P}{\cH}{Q} \cong \rL^2(P) \ot_{P_2} \cH'$, where $P_2 = A \rtimes_\Omega G_1$ and where the $P_2$-$Q$-bimodule $\cH'$ is defined as
\begin{equation}\label{eq.whatweget}
\begin{split}
& \bim{P_2}{\cH'}{Q} = \bim{\theta(P_2)}{(q (\C^n \ot \rL^2(Q)) \ot K)}{Q} \quad\text{with}\\ & \theta : P_2 \recht q(\M_n(\C) \ot Q)q \ot \B(K) : \theta(au_g) = \psi(au_g) \ot \pi(g) \; .
\end{split}
\end{equation}
Take the non-negligible subset $\cU \subset \{1,\ldots,n\} \times X$ such that $q = \chi_\cU$ and take the isomorphism $\Delta : X \recht \cU$ such that $\psi(a) = a \circ \Delta^{-1}$ for all $a \in A$. Denote by $\Delta_1 : X \recht X$ the composition of $\Delta$ and $(i,x) \mapsto x$. It follows that $\Delta_1$ is a local isomorphism, locally multiplying the measure by $(\Tr \ot \tau)(q)$ and satisfying $\Delta_1(G_1 \cdot x) \subset H \cdot \Delta_1(x)$ almost everywhere. By Lemma \ref{lem.quotient}, we find an $m$-to-$1$ quotient map $\Delta_2 : X \recht X$ and a group homomorphism $\delta : G_1 \recht H$ such that $\Delta_2(g \cdot x) = \delta(g) \cdot \Delta_2(x)$ and $\Delta_2(x) \in H \cdot \Delta_1(x)$ almost everywhere. Moreover, Lemma \ref{lem.quotient} provides a subgroup $\Lambda < \Aut(X,\mu)$ of order $m$, with $g \Lambda g^{-1} = \Lambda$ for all $g \in G_1$ and such that $\Delta_2$ induces a conjugacy between the actions
$$\frac{G_1}{G_1 \cap \Lambda} \actson \frac{X}{\Lambda} \quad\text{and}\quad H_1 \actson X \quad\text{where}\quad H_1 = \delta(G_1) \; .$$
Since $\Delta_2(x) \in H \cdot \Delta_1(x)$ and since $\Delta_2$ is an $m$-to-$1$ quotient map, it follows that $(\Tr \ot \tau)(q) = m$ and that there exists a $W \in q(\M_{n,m}(\C) \ot Q)$ satisfying $WW^* = q$, $W^* W = 1$ and, writing $\psi_1(\,\cdot \,) = W^* \psi(\,\cdot\,)W$, $\psi_1(a) = 1 \ot a \circ \Delta_2^{-1}$ for all $a \in \rL^\infty(X/\Lambda) \subset A$. Since $\Delta_2(g \cdot x) = \delta(g) \cdot \Delta_2(x)$, it then follows that $\psi_1(u_g) \in (\M_m(\C) \ot A) u_{\delta(g)}$ for all $g \in G_1$. In particular, $\psi_1(P_1) \subset \M_m(\C) \ot (A \rtimes_{\omega} H_1)$. Since $\psi_1(P_1) \subset q(\M_n(\C) \ot Q)q$ has finite index, this implies that $H_1 < H$ has finite index.

Since $H_1 < H$ has finite index, the action $H_1 \actson (X,\mu)$ is cocycle superrigid. But $H_1 \actson X$ is conjugate with $G_1/(G_1 \cap \Lambda) \actson X/\Lambda$. Hence, Lemma \ref{lem.noquotient} implies that $m=1$, $\Lambda = \{e\}$. This means that $\Delta_2 \in \Aut(X,\mu)$ and $\Delta_2 g \Delta_2^{-1} = \delta(g)$ for all $g \in G_1$. Hence, $\Delta_2 \in \cG$. Put $\gamma = \Delta_2^{-1}$. So, we may from now on assume that $n=1$, $q=1$ and that $\psi_1 : P_1 \recht Q$ is given by
$$\psi_1(a) = a(\gamma \cdot \,) \quad\text{and}\quad \psi_1(u_g) = a_{\gamma^{-1} g \gamma} u_{\gamma^{-1} g \gamma}$$
for some unitaries $(a_h)_{h \in H_1} \in \cU(A)$.

Define the measurable map $\eta : H_1 \times X \recht \T : \eta(h,x) = a_h(h \cdot x)$. It follows that
$$\eta(hk,x) = \bigl(\omega\,\overline{\Omega_1 \circ \Ad \gamma} \bigr)(h,k) \, \eta(h, k \cdot x) \, \eta(k,x) \; .$$
Lemma \ref{lem.untwist-pseudo} yields a unitary $b \in \cU(A)$ and a map $\rho : G_1 \recht \T$ such that, replacing $\psi_1$ by $\Ad b \circ \psi_1$, we have $\psi_1(a) = a(\gamma \cdot \,)$ for all $a \in A$ and $\psi_1(u_g) = \rho(g) u_{\gamma^{-1} g \gamma}$ for all $g \in G_1$.

Replacing the projective representation $\pi : G_1 \recht \cU(K)$ by the projective representation $G_1 \recht \cU(K) : g \mapsto \rho(g) \pi(g)$, we may assume that $\psi(u_g) = u_{\gamma^{-1} g \gamma}$ for all $g \in G_1$. In particular, $\Omega_\pi = \Omega \overline{\om \circ \Ad \gamma^{-1}}$.

Put $G_0 = G \cap \gamma H \gamma^{-1}$ and $H_0 = H \cap \gamma^{-1} G \gamma$. The $2$-cocycle $\Omega \overline{\om \circ \Ad \gamma^{-1}}$ makes sense on $G_0$ and hence is the obstruction $2$-cocycle of an induction of $\pi$ to a projective representation $\pitil$ of $G_0$ (see e.g.\ \cite[Definition 6.8]{Va07}). Put $P_0 = A \rtimes_{\Omega} G_0$ and $Q_0 = A \rtimes_\omega H_0$. From \eqref{eq.whatweget} and the isomorphism $\cH \cong \rL^2(P) \ot_{P_1} \cH'$, one deduces that
$$\bim{P}{\cH}{Q} \cong \rL^2(P) \ot_{P_0} \cK(\gamma,\pitil) \ot_{Q_0} \rL^2(Q) \; ,$$
where $\cK(\gamma,\pitil)$ is as in the formulation of Theorem \ref{thm.main}.2. Since $\cH$ was assumed to be irreducible, it follows that $\pitil$ is irreducible.
\end{proof}

\begin{proof}[{\bf Proof of Theorem \ref{thm.main}, point 1.}]
Put $M = A \rtimes \Gamma$. Let $Q$ be a II$_1$ factor and $\bim{M}{\cK}{Q}$ a non-zero finite index bimodule.
Define $\cH = \cK \ot_Q \overline{\cK}$. Denote by $\cG$ the commensurator of $\Gamma$ inside $\Aut(X,\mu)$. Since $\cH$ is a non-zero finite index $M$-$M$-bimodule, the already proven point 2 of Theorem \ref{thm.main} yields a finite dimensional Hilbert space $K$ and a finite index inclusion $\psi_0 : M \recht \B(K) \ot M$ such that $\bim{M}{\cH}{M} \cong \bim{\psi_0(M)}{(K \ot \rL^2(M))}{M}$ and such that the restriction $\psi_0|_A$ has the special form
\begin{equation}\label{eq.specialform}
\psi_0(a) = \sum_{i=1}^n p_i \ot a(\gamma_i \, \cdot) \quad\text{where}\;\; \gamma_1,\ldots,\gamma_n \in \cG \;\;\text{and}\;\; p_1,\ldots,p_n \;\;\text{are projections with sum $1$.}
\end{equation}
After a unitary conjugacy and a regrouping of the $p_i$ and $\gamma_i$, we may assume that $\gamma_1 \Gamma,\ldots,\gamma_n \Gamma$ are mutually disjoint. By construction, we find an intermediate subfactor $\psi_0(M) \subset P \subset \B(K) \ot M$ such that $P$ and $Q$ are stably isomorphic. For the rest of the proof, we only retain the information that $Q$ is stably isomorphic with an intermediate subfactor of a finite index inclusion $\psi_0 : M \recht \B(K) \ot M$ where $\psi_0|_A$ is of the special form \eqref{eq.specialform}.

Put $D = \psi_0(A)' \cap (\B(K) \ot M)$. Since $\cG$ acts freely on $(X,\mu)$, we have $D = \bigoplus_i \B(p_iK) \ot A$.
Denote by $\Gamma_0 < \Gamma$ a finite index subgroup with the property that $\gamma_i^{-1} \Gamma_0 \gamma_i \subset \Gamma$ for all $i=1,\ldots,n$. For every $g \in \Gamma_0$, define
$$v_g := \sum_{i=1}^n p_i \ot u_{\gamma_i^{-1} g \gamma_i} \; .$$
A direct computation shows that $\psi_0(u_g) v_g^*$ commutes with $\psi(A)$ and hence, $\psi_0(u_g) \in D v_g$ for all $g \in \Gamma_0$.
Write $\psi_0(u_g) = \sum_{i=1}^n \eta_{i,g}(p_i \ot u_{\gamma_i^{-1} g \gamma_i})$ for some unitaries $\eta_{i,g} \in \B(p_iK) \ot A$. For a fixed $i$, these unitaries can be reinterpreted as a $1$-cocycle for the action $\Gamma_0 \actson X$ with values in $\cU(p_i K)$. As a result, we can unitarily conjugate $\psi_0$ and assume that $\psi_0(u_g) = \sum_{i=1}^n \pi_i(g)p_i \ot u_{\gamma_i^{-1} g \gamma_i}$ for all $g \in \Gamma_0$. Since the action of $\gamma_i^{-1} \Gamma_0 \gamma_i$ on $(X,\mu)$ is weakly mixing, it follows that
$$\psi_0(M)' \cap (\B(K) \ot M) \subset \bigoplus_i (\B(p_i K) \ot 1) \; .$$
We replace $\psi_0(M)$ by $\psi_0(M)q$ where $q$ is a minimal projection in $P \cap \psi_0(M)'$. So, $\psi_0|_A$ still has the special form \eqref{eq.specialform} and also the $\psi_0(u_g), g \in \Gamma_0$ keep their special form. We now have that $\psi_0(M) \subset P$ is irreducible.

We apply Lemma \ref{lem.projrep} to the irreducible inclusion $\psi_0(M) \subset P$. We find a finite index subgroup $\Gamma_1 < \Gamma$, a projection $p \in P \cap \psi(A)'$, an obstruction $2$-cocycle $\Omega \in \rZ^2(\Gamma_1,\T)$ of a finite dimensional projective representation and a finite index inclusion $\psi : A \rtimes_\Omega \Gamma_1 \recht pPp$ such that $\psi(a) = \psi_0(a) p$ for all $a \in A$ and $\psi(A) \subset pPp$ is maximal abelian. Moreover, we can take $p$ such that $\lspan \{\psi_0(u_g) p \psi_0(u_g)^* \mid g \in \Gamma\}$ is finite dimensional. Since $\Gamma \actson (X,\mu)$ is weakly mixing, we conclude that $p \in \bigoplus_i (\B(p_i K) \ot 1)$.

Replace $P$ by $pPp$ and $K$ by $pK$. We have found a finite index inclusion $\psi : A \rtimes_\Omega \Gamma_1 \recht \B(K) \ot M$ such that $P$ is an intermediate subfactor $\psi(A \rtimes_\Omega \Gamma_1) \subset P \subset \B(K) \ot M$ with $\psi(A) \subset P$ being maximal abelian and such that $\psi|_A$ is of the special form \eqref{eq.specialform}.

Define the group
$$G = \{\gamma \in \cG \mid \;\text{there exists a unitary $w \in \cU(P)$ such that}\;\; w^* \psi(a) w = \psi(a(\gamma \, \cdot))\;\;\forall a \in A \} \; .$$
By construction $\Gamma_1 < G$. Also, since $\cG$ acts freely on $(X,\mu)$, if $\gamma \in G - \Gamma_1$ and $w \in \cU(P)$ satisfies $w^* \psi(a) w = \psi(a(\gamma\, \cdot))$ for all $a \in A$, then $E_{\psi(A \rtimes_\Omega \Gamma_1)}(w) = 0$. Because $\psi(A \rtimes_\Omega \Gamma_1) \subset P$ has finite index, it follows that $\Gamma_1 < G$ has finite index. Hence, $G$ and $\Gamma$ are commensurate groups inside $\Aut(X,\mu)$. Making $\Gamma_1$ smaller if necessary, we may assume that $\Gamma_1$ is a normal subgroup of $G$.

Let $\gamma \in G$. Take a unitary $w \in \cU(P)$ satisfying $w^* \psi(a) w = \psi(a(\gamma\, \cdot))$ for all $a \in A$. As such, $w$ is determined up to multiplication by a unitary in $\psi(A)$. We claim that $w$ can be chosen in such a way that $w^* \psi(u_g) w \in  \T \psi(u_{\gamma^{-1} g \gamma})$ for all $g \in \Gamma_1$. To prove this claim, first observe that $w^* \psi(u_g) w \psi(u_{\gamma^{-1} g \gamma})^*$ commutes with $\psi(A)$ and hence, belongs to $\psi(A)$. We find unitaries $(a_g)_{g \in \Gamma_1}$ in $A$ such that
$$w^* \psi(u_g) w = \psi(a_{\gamma^{-1} g \gamma} \, u_{\gamma^{-1} g \gamma}) \quad\text{for all}\;\; g \in \Gamma_1 \; .$$
Lemma \ref{lem.untwist-pseudo} yields a unitary $b \in \cU(A)$ such that, replacing $w$ by $w \psi(b)$, the claim is proven.

By weak mixing of $\Gamma_1 \actson (X,\mu)$, the unitary $w$ satisfying the claim in the previous paragraph, is uniquely determined up to multiplication by a scalar. So, for every $g \in G$, choose $w_g$ satisfying the claim. Make this choice such that $w_g = \psi(u_g)$ for all $g \in \Gamma_1$. By uniqueness of the $w$, we get $w_g w_h = \Omega(g,h) w_{gh}$ for all $g,h \in G$ and for some $\Omega \in \rZ^2(G,\T)$ that extends the given $\Omega \in \rZ^2(\Gamma_1,\T)$. Since $E_{\psi(A)}(w_g) = 0$ when $g \neq e$, the formula
$$\theta : A \rtimes_\Omega G \recht P : \theta(a) = \psi(a) \;\;\text{for}\;\; a \in A \quad\text{and}\quad \psi(u_g) = w_g \;\;\text{for}\;\;g \in G \; ,$$
is a well defined finite index inclusion. Since $\Omega|_{\Gamma_1}$ was the obstruction $2$-cocycle of a finite dimensional projective representation, the same is true for $\Omega \in \rZ^2(G,\T)$ by considering the induced projective representation.

It remains to prove that $\theta$ is surjective. Denote $P_0 = \theta(A \rtimes_\Omega G)$, which equals the von Neumann subalgebra of $P$ generated by $\psi(A)$ and the unitaries $w_g, g \in G$. Choose $i,j$, an operator $T \in p_i \B(K) p_j$, $b \in A$, $g \in \Gamma$. Define $d = E_P(T \ot b u_g)$. It suffices to prove that $d \in P_0$. We may assume that $d \neq 0$. Denote $\gamma = \gamma_i g \gamma_j^{-1}$. It follows that $\psi(a) d = d \psi(a(\gamma \, \cdot))$ for all $a \in A$. Let $d = v |d|$ be the polar decomposition of $d$. Then, $|d| \in \psi(A) \subset P_0$, while $v$ is a non-zero partial isometry satisfying $\psi(a) v = v \psi(a(\gamma \, \cdot))$ for all $a \in A$. Let $\cU \subset X$ be the non-negligible subset such that $v^* v = \psi(\chi_\cU)$. Put $\Gamma_2 = \Gamma_1 \cap \gamma^{-1} \Gamma_1 \gamma$. Since $\Gamma_2$ acts ergodically on $(X,\mu)$, we can find subsets $\cU_n \subset \cU$ and group elements $g_n \in \Gamma_2$ such that the sets $g_n^{-1} \cdot \cU_n$ form a partition of $X$, up to measure zero. We may assume that $\cU_0 = \cU$ and $g_0 = e$. It is then easy to check that
$$w := \sum_n \psi(u_{\gamma g_n \gamma^{-1}})^* \; v \; \psi(\chi_{\cU_n} \, u_{g_n})$$
is a unitary in $P$ satisfying $\psi(a) w = w \psi(a(\gamma \, \cdot))$ for all $a \in A$. It follows that $\gamma \in G$ and $w \in P_0$. Since $v = w \psi(\chi_\cU)$, we also get $v \in P_0$. Hence, $d \in P_0$, ending the proof of the theorem.
\end{proof}

\section{Concrete computations of the index sets $\cC(M)$ and $\cI(M)$}

\begin{theorem} \label{thm.concrete}
Consider $\Q \subset K$ where $K$ is a countable field of characteristic zero and $\Q \neq K$. Define $\Gamma_1 = \SL(3,\Q)$ and $\Gamma_2 = \SL(3,K)$. Let $q \in \Q \setminus \{0,1,-1\}$ and define
$$\Sigma := A^\Z \quad\text{where}\quad A = \begin{pmatrix} 1 & 0 & 0 \\ 1 & q & 0 \\ 1 & 0 & q^{-1} \end{pmatrix} \; .$$
We consider $\Sigma$ as a common subgroup of $\Gamma_1$, $\Gamma_2$ and put $\Gamma = \Gamma_1 *_\Sigma \Gamma_2$. Finally, consider the subgroup $\Lambda := \Lambda_1 * \Lambda_2 < \Gamma$, where both $\Lambda_i$ are given by
$$\left\{ \begin{pmatrix} 1 & a & b \\ 0 & 1 & c \\ 0 & 0 & 1 \end{pmatrix} \; \Big| \; a,b,c \in \Z \right\} \; .$$
For any atomic probability space $(X_0,\mu_0)$ with unequal weights $\mu_0(\{x\})$, we consider the II$_1$ factor
$$M := \rL^\infty\bigl( X_0^{\Gamma/\Lambda} \bigr) \rtimes \Gamma \; .$$
Then, $\cC(M) = \{|G| \mid G < \Gal(K) \;\;\text{is a finite subgroup}\;\}$.
\end{theorem}

\begin{corollary} \label{cor.concrete}
Let $\cN \subset \N$ be a subset of the natural numbers that is closed under taking divisors: if $n \in \cN$ and $m | n$, then $m \in \cN$.
There exists a countable field $K$ of characteristic zero such that the associated II$_1$ factor $M$ constructed in Theorem \ref{thm.concrete} satisfies $\cC(M) = \cN$.
\end{corollary}

\begin{proof}[Proof of Theorem \ref{thm.concrete}]
We denote $(X,\mu) := (X_0,\mu_0)^{\Gamma/\Lambda}$. Note that $\Gamma$ has no non-trivial finite dimensional unitary representations. In particular, $\Gamma$ has no non-trivial finite index subgroups.

Throughout the proof, we use the following notation: whenever $H_1,H_2 < G$ are subgroups, we write $H_1 \prec_G H_2$ if there exists a $g \in G$ such that $H_1 \cap g^{-1} H_2 g$ has finite index in $H_1$. In other words, a finite index subgroup of $H_1$ can be conjugated into $H_2$.

Also note that whenever $H < \Gamma_i$ is a subgroup such that $H \not\prec_{\Gamma_i} \Sigma$ and $H \not\prec_{\Gamma_i} \Lambda_i$, then $H$ acts with infinite orbits on $\Gamma/\Lambda$ and hence, $H$ acts weakly mixingly on $(X,\mu)$. This applies to $H = \SL(3,\Z)$ and its finite index subgroups. It also applies to $H = E_{ij}(\Q)$, where $E_{ij}(x) = 1 + e_{ij}(x)$ with $e_{ij}(x)$ being the obvious elementary matrix with a single non-zero entry given by $x$.

We call reduced expression every product of elements alternatingly from $\Gamma_1 - \Sigma$ and $\Gamma_2 - \Sigma$. Every $g \in \Gamma - \Sigma$ admits a reduced expression. We refer to the first factor of such an expression as the first letter of $g$. The first letter of $g$ is uniquely determined up to right multiplication by an element from $\Sigma$. The length of a reduced expression for $g$ is denoted by $|g|$. By convention $|g|=0$ when $g \in \Sigma$.

We also have that $\Sigma$ acts with infinite orbits on $\Gamma/\Lambda$. If not, we would find a finite index subgroup $\Sigma_0 < \Sigma$ and a $g \in \Gamma$ such that $g \Sigma_0 g^{-1} \subset \Lambda$. Since $\Sigma \cap \Lambda = \{e\}$, we have $g \not\in \Lambda \Sigma$. Hence, we can write $g = g_0 g_1 \cdots g_n$ with $g_0 \in \Lambda$ and $g_1 \cdots g_n$ being a reduced expression with $g_1 \in \Gamma_i - \Lambda_i \Sigma$. Since $g_1 \cdots g_n$ conjugates $\Sigma_0$ into $\Lambda$ and $g_1 \in \Gamma_i - \Lambda_i \Sigma$, we conclude that $(g_k \cdots g_n) \Sigma_0 (g_k \cdots g_n)^{-1} \subset \Sigma$, first for $k=n$, then for $k=n-1$, until $k=2$. Let $h$ be a non-trivial element in $(g_2 \cdots g_n) \Sigma_0 (g_2 \cdots g_n)^{-1}$. Then, $h \in \Sigma$ and $g_1 h g_1^{-1} \in \Lambda$. So, $g_1 h g_1^{-1} \in \Lambda_i$. The spectrum of every element in $\Lambda_i$ is $\{1\}$ and hence, the spectrum of $h$ equals $\{1\}$. This is a contradiction with $h \in \Sigma - \{e\}$.

{\bf Part 1.} The action $\Gamma \actson (X,\mu) = (X_0,\mu_0)^{\Gamma/\Lambda}$ is cocycle superrigid for arbitrary countable or compact second countable target groups. This is, as follows, a direct consequence of results in \cite{Po05}. Assume that $\om : \Gamma \times X \recht \cG$ is a $1$-cocycle with values in the discrete or compact second countable group $\cG$. We first claim that both restrictions $\om|_{\Gamma_i}$, $i=1,2$, are cohomologous to a group morphism. We prove this claim for $i=2$, the case $i=1$ being analogous. Since $\SL(3,\Z)$ has property (T), \cite[Theorem 0.1]{Po05} allows us to assume that $\om|_{\SL(3,\Z)}$ is a group morphism. Since $\SL(3,\Z)$ acts with infinite orbits on $\Gamma/\Lambda$, \cite[Proposition 3.6]{Po05} implies that $\om$ is a group morphism on $\SL(3,\Q)$. Since $E_{ij}(\Q)$ acts with infinite orbits on $\Gamma/\Lambda$, the same \cite[Proposition 3.6]{Po05} implies that $\om$ is a group morphism on $E_{ij}(K)$. Since the $E_{ij}(K)$ generate $\Gamma_2$, it follows that $\om$ is a group morphism on $\Gamma_2$.

So, we have proven the claim and may assume that $\om(g,x) = \delta_1(g)$ when $g \in \Gamma_1$, while $\om(g,x) = \vphi(g \cdot x) \delta_2(g) \vphi(x)^{-1}$ when $g \in \Gamma_2$. Since $\Sigma$ acts with infinite orbits on $\Gamma/\Lambda$, it follows that $\vphi$ is essentially constant (see e.g.\ \cite[Lemma 5.4]{PV08b}) and hence, $\om$ is a group morphism, concluding the proof of part 1.

{\bf Part 2.} Every finite index $M$-$M$-bimodule contains a non-zero finite index $\rL^\infty(X)$-$\rL^\infty(X)$-subbimodule. This is a consequence of \cite[Theorem 7.1]{PV09}. Denote by $D_i$ the normalizer of $\Sigma$ inside $\Gamma_i$. We start with the following observation: whenever $g \in \Gamma_i$ such that $g \Sigma g^{-1} \cap \Sigma \neq \{1\}$, we have $g \in D_i$. Indeed, if $h \in \Sigma - \{1\}$ and $g h g^{-1} \in \Sigma$, consider the spectrum of $h$ and conclude that $g h g^{-1}$ equals either $h$ or $h^{-1}$. This moreover implies that $g A g^{-1}$ equals either $A$ or $A^{-1}$ and hence $g \in D_i$. The observation implies that $g \Sigma g^{-1} \cap \Sigma = \{1\}$ whenever $g \in \Gamma - D_1 *_\Sigma D_2$. So, \cite[Theorem 7.1]{PV09} yields the result.

{\bf Part 3.} The normalizer of $\Gamma$ inside $\Aut(X,\mu)$ is generated by $\Gamma$ and the group $\Aut(\Lambda < \Gamma) := \{\delta \in \Aut(\Gamma) \mid \delta(\Lambda) = \Lambda \}$ which acts on $(X,\mu)$ by $(\delta \cdot x)_{g\Lambda} = x_{\delta^{-1}(g) \Lambda}$. Part 3 follows from \cite[Lemma 6.15]{PV09}, once we have shown that $\Lambda$ acts with infinite orbits on $\Gamma/\Lambda - \{e \Lambda\}$. It is easy to check that $\Lambda g \Lambda$ is an infinite subset of $\Gamma/\Lambda$ whenever $g \in \Gamma - \Sigma \Lambda$. If $\si \in \Sigma - \{e\}$ and if $\Lambda \si \Lambda$ would be a finite subset of $\Gamma/\Lambda$, we find a finite index subgroup $\Lambda_0 < \Lambda$ such that $\si^{-1} \Lambda_0 \si \subset \Lambda$. In particular, $\si^{-1}(\Lambda_0 \cap \Lambda_1) \si \subset \Lambda_1$. Hence, $\si$ is upper triangular (with arbitrary diagonal elements). So, the standard basis vector $e_1$ is an eigenvector of $\si$. But all non-trivial elements of $\Sigma$ have the same set of eigenvectors, namely the non-zero multiples of the vectors $(1-q)e_1 + e_2 - qe_3$, $e_2$ and $e_3$. The vector $e_1$ is not in this set, yielding the required contradiction.

{\bf Part 4.} Every automorphism $\delta \in \Aut(\Gamma)$ is of the form $\delta = (\Ad g) \circ (\delta_1 * \delta_2)$ where $g \in \Gamma$ and $\delta_i \in \Aut(\Gamma_i)$ satisfy $\delta_i(\Sigma) = \Sigma$ and $(\delta_1)|_\Sigma = (\delta_2)|_\Sigma$.

We start the proof of part 4 with the following preliminary statement on general amalgamated free products $\Gamma = \Gamma_1 *_\Sigma \Gamma_2$. Whenever $G < \Gamma$ is a subgroup, denote $|G| := \sup \{|g| \mid g \in G \}$. One can easily show that $|G| < \infty$ if and only if there exists $g \in \Gamma$ and $i \in \{1,2\}$ such that $g G g^{-1} \subset \Gamma_i$. One direction being obvious, assume that $G < \Gamma$ is a subgroup and $|G| < \infty$. Take $g \in \Gamma$ minimizing the function $g \mapsto |gGg^{-1}|$. Replace $G$ by $g G g^{-1}$. We show that $G < \Gamma_i$ for some $i \in \{1,2\}$. If $|G| \leq 1$, we get $G \subset \Gamma_1 \cup \Gamma_2$. Since $G$ is a subgroup, one easily checks that $G$ actually sits in one of the $\Gamma_i$. So, assume that $|G| \geq 2$. We will produce an element $h_1$ such that $|h_1^{-1} G h_1| < |G|$, which contradicts our minimal choice of $g$. Take $h \in G$ with $|h| = |G|$. Write $h$ in reduced form and denote by $h_1$ the first letter of $h$. Assume that $h_1 \in \Gamma_1 - \Sigma$. We claim that every element $k \in G$ either belongs to $\Gamma_1$ or admits a reduced expression starting with $h_1$ and ending with $h_1^{-1}$. Let $k \in G$ and $k \not\in \Gamma_1$. If the first letter of $k$ cannot be chosen to be $h_1$, then $|k^{-1} h | > |h| = |G|$, which is absurd. If the last letter of $k$ cannot by chosen to be $h_1^{-1}$, we have $|k h| > |G|$. This proves the claim and hence, $|h_1^{-1} G h_1| < |G|$.

We now return to the actual proof of part 4. Whenever $0 < \rho < 1$, the function $\vphi_\rho(g) = \rho^{|g|}$ is positive definite. If $\rho \recht 1$, the functions $\vphi_\rho$ tend to $1$ pointwise. Consider $\SL(3,\Z)$ as a subgroup of $\Gamma_2$. By property (T), $\vphi_\rho$ converges uniformly to $1$ on $\delta(\SL(3,\Z))$. This means that $|\delta(\SL(3,\Z))| < \infty$. So, after replacing $\delta$ by $\Ad g \circ \delta$, we may assume that $\delta(\SL(3,\Z))$ is a subgroup of either $\Gamma_1$ or $\Gamma_2$. We assume that $\delta(\SL(3,\Z)) < \Gamma_2$ and explain later why the other option is impossible.

Observe that whenever $G < \Gamma_i$ such that $G \not\prec_{\Gamma_i} \Sigma$ and whenever $g \in \Gamma$ satisfies $g G g^{-1} \subset \Gamma_i$, then $g \in \Gamma_i$.

Whenever $g \in \SL(3,\Q)$ viewed as a subgroup of $\Gamma_2$, $g$ quasi-normalizes $\SL(3,\Z)$. So, $\delta(g)$ quasi-normalizes the subgroup $\delta(\SL(3,\Z))$ of $\Gamma_2$. By the observation in the previous paragraph, $\delta(g) \in \Gamma_2$. So, $\delta(\SL(3,\Q)) \subset \Gamma_2$. For all $i \neq j$, we have that $\delta(E_{ij}(K))$ commutes with $\delta(E_{ij}(\Q)) < \Gamma_2$. Again applying the observation in the previous paragraph, we get that $\delta(E_{ij}(K)) \subset \Gamma_2$. We have shown that $\delta(\Gamma_2) \subset \Gamma_2$.

Similarly, we find $h \in \Gamma$ and $i \in \{1,2\}$ such that $\delta(\Gamma_1) \subset h \Gamma_i h^{-1}$. Since $\delta(\Gamma_1)$ and $\delta(\Gamma_2)$ together generate $\Gamma$, we get that $i = 1$ and $h \in \Gamma_2 \Gamma_1$. Replacing $\delta$ by $\Ad h_0 \circ \delta$ for some $h_0 \in \Gamma_2$, we have found that $\delta(\Gamma_i) < \Gamma_i$ for both $i=1,2$. Since $\delta$ is surjective, we finally find that $\delta(\Gamma_i) = \Gamma_i$ for both $i = 1,2$. It automatically follows that $\delta_i(\Sigma) = \Sigma$ and $(\delta_1)|_\Sigma = (\delta_2)|_\Sigma$.

If in the beginning $\delta(\SL(3,\Z))$ would belong to a conjugate of $\Gamma_1$, the argument would go through in exactly the same way and end up with finding an isomorphism between $\SL(3,K)$ and $\SL(3,\Q)$. Since $\Q \neq K$, this is impossible.

{\bf Part 5.} We have $\Aut(\Lambda < \Gamma) = \Ad \Lambda \times \Gal(K)$, where $\al \in \Gal(K)$ defines the automorphism $\theta_\al \in \Aut(\Gamma)$ that is the identity on $\Gamma_1$ and the pointwise application of $\al$ on $\Gamma_2$.

Let $\delta \in \Aut(\Lambda < \Gamma)$. By part 4, take $g \in \Gamma$ and $\delta_i \in \Aut(\Gamma_i)$ such that $\delta_i(\Sigma) = \Sigma$, $(\delta_1)|_\Sigma = (\delta_2)|_\Sigma$ and $\delta = \Ad g \circ (\delta_1 * \delta_2)$. Since $\Lambda_i \not\prec_{\Gamma_i} \Sigma$ and $\delta_i(\Sigma) = \Sigma$, also $\delta(\Lambda_i) \not\prec_{\Gamma_i} \Sigma$. Hence, the fact that $g \delta_i(\Lambda_i) g^{-1} \subset \Lambda$, forces $g \in \Lambda \Gamma_i$ for both $i = 1,2$. So, $g \in \Lambda \Sigma$. Making the appropriate replacements, we may assume that $\delta_i(\Sigma) = \Sigma$ and $\delta_i(\Lambda_i) = \Lambda_i$. It remains to prove that $\delta_1 = \id$ and $\delta_2 = \theta_\al$ for some $\al \in \Gal(K)$.

Denote by $\beta$ the automorphism of $\SL(3,K)$ given by the composition of the inverse and the transpose: $\beta(g) = (g^T)^{-1}$.
As an automorphism of $\SL(3,K)$, we either have $\delta_2 = \Ad B^{-1} \circ \theta_\al$ or we have $\delta_2 = \Ad B^{-1} \circ \theta_\al \circ \beta$ for some $B \in \GL(3,K)$ and $\al \in \Gal(K)$. In the latter case, it would follow that $B^{-1} \beta(\Sigma) B = \Sigma$ and $B^{-1} \beta(\Lambda_2) B = \Lambda_2$. The formula $B^{-1} \beta(\Lambda_2) B = \Lambda_2$ implies that $B$ is of the form
$\bigl(\begin{smallmatrix} 0 & 0 & 1 \\ 0 & 1 & 0 \\ 1 & 0 & 0 \end{smallmatrix} \bigr)\bigl(\begin{smallmatrix} * & * & * \\ 0 & * & * \\ 0 & 0 & * \end{smallmatrix} \bigr) = \bigl(\begin{smallmatrix} 0 & 0 & * \\  0 & * & * \\ * & * & *  \end{smallmatrix} \bigr)$. We only retain that $B_{12} = 0$.

The non-trivial elements of $\Sigma$ all have the same eigenvectors, namely the non-zero multiples of the vectors $(1-q)e_1 + e_2 - qe_3$ (with eigenvalue $1$), and $e_2$, $e_3$ (with eigenvalue a non-zero power of $q$). On the other hand, the non-trivial elements of $\Sigma^T$ have as eigenvectors the non-zero multiples of $e_1$ (with eigenvalue $1$) and $e_1 + (q-1)e_2, qe_1 + (1-q)e_3$ (with eigenvalue a non-zero power of $q$). Whenever $\sigma \in \Sigma - \{e\}$, the formula $B^{-1} \beta(\Sigma) B = \Sigma$ implies that $B$ maps the eigenvectors of $\sigma$ to eigenvectors of $\beta(\sigma) \in \Sigma^T$ with the same eigenvalue. Hence, $Be_2$ must be a non-zero multiple of either $e_1 + (q-1)e_2$ or $qe_1 + (1-q)e_3$. In both cases $B_{12} \neq 0$, yielding a contradiction.

Hence, $\delta_2 = \Ad B^{-1} \circ \theta_\al$ for some $B \in \GL(3,K)$ and $\al \in \Gal(K)$. So, $B^{-1} \Lambda_2 B = \Lambda_2$, implying that $B$ is upper triangular. Also, $B^{-1} \Sigma B = \Sigma$, implying that $B^{-1} A B = A$ or $B^{-1} A B = A^{-1}$. In the latter case $B$ must map $e_2$ to a non-zero multiple of $e_3$, contradicting the fact that $B$ is upper triangular. So, $A$ and $B$ commute. It follows that the eigenvectors $(1-q)e_1 + e_2 - qe_3$, $e_2$, $e_3$ of $A$ are all eigenvectors of $B$. Since $B$ is upper triangular, also $e_1$ is an eigenvector of $B$. Both together imply that $B$ is a scalar multiple of the identity matrix and hence, $\delta_2 = \theta_\al$.

The proof that $\delta_1 = \id$ is identical to the proof that $\delta_2 = \theta_\al$.

{\bf Conclusion.} The group $\Gamma$ has no non-trivial finite index subgroups, so that part 1 and 2 say that $\Gamma \actson (X_0,\mu_0)^{\Gamma/\Lambda}$ satisfies Condition \ref{cond}. By parts 3 and 5, the normalizer of $\Gamma$ inside $\Aut(X,\mu)$ is given by $\Gamma \rtimes \Gal(K)$. So, Corollary \ref{cor.case-no-fin-dim-rep} provides the required formula for $\cC(M)$.
\end{proof}

\begin{proof}[Proof of Corollary \ref{cor.concrete}]
We mimic the construction of \cite[Theorem 4.7.1]{We06}.
Whenever $p$ and $q$ are prime numbers with $p | q-1$, denote by $\xi_q$ a primitive $q$-th root of unity and choose a subgroup $G$ of order $(q-1)/p$ in $\Gal(\Q(\xi_q)) \cong \mathbb{F}_q^\times \cong \frac{\Z}{(q-1)\Z}$. Define $K_{p,q}$ to be the fixed subfield $\Q(\xi_q)^G$ of $\Q(\xi_q)$. By the fundamental theorem of Galois theory, $\Gal(K_{p,q}) \cong \frac{\Z}{p\Z}$.

Denote by $\cP$ the set of prime numbers. For every $p \in \cP$, put $f(p) = \sup\{k \in \N \mid p^k \in \cN\}$. By convention, we put $f(p) = 0$ if $p \not\in \cN$ and we put $f(p) = +\infty$ if $p^k \in \cN$ for all $k$.

Let $p_n$ be a finite or infinite sequence of prime numbers with every prime number $p$ appearing $f(p)$ times. Observe that $\cN$ consists of those natural numbers that divide $p_0 \cdots p_n$ for $n$ large enough.

We can inductively choose distinct prime numbers $q_n$ such that $q_n = 1 \mod p_n$. Define $K$ as the subfield of $\C$ generated by all the $K_{p_n,q_n}$. As e.g.\ in \cite[Theorem 4.7.1]{We06}, one easily checks that $\Gal(K)$ is isomorphic with the compact group $\prod_n \frac{\Z}{p_n \Z}$.

Define $G_n = \prod_{k \leq n} \frac{\Z}{p_k \Z}$ and denote by $\pi_n : \Gal(K) \recht G_n$ the natural quotient map. Whenever $G < \Gal(K)$ is a finite subgroup, there exists an $n$ such that $\pi_n$ is injective on $G$. Hence, $|G|$ divides $p_0 \cdots p_n$. Conversely, if $m$ divides $p_0 \cdots p_n$, one easily defines a subgroup $G < \Gal(K)$ with $|G| = m$.
\end{proof}

\section{The C$^*$-tensor category $\Bimod(M)$ and the C$^*$-bicategory of commensurable II$_1$ factors}\label{sec.categorical}

Whenever $G \actson (X,\mu)$ satisfies Condition \ref{cond}, Theorem \ref{thm.main} provides a description of all II$_1$ factors that are commensurable with $\rL^\infty(X) \rtimes G$ and moreover gives a description of all finite index bimodules between these commensurable II$_1$ factors.

Given a II$_1$ factor $M$, we denote by $\Comm(M)$ the C$^*$-bicategory of commensurable II$_1$ factors. The objects in this category (also called $0$-cells) are all II$_1$ factors that are commensurable with $M$. The morphisms ($1$-cells) are given by the finite index bimodules and their composition is defined by the Connes tensor product. The morphisms between the morphisms ($2$-cells) are given by the bimodular bounded Hilbert space operators. The bicategory $\Comm(M)$ in particular encodes the C$^*$-tensor category $\Bimod(M)$ of all finite index $M$-$M$-bimodules equipped with the Connes tensor product.

The aim of this section is to reinterpret Theorem \ref{thm.main} and to provide an explicit description of the C$^*$-bicategory $\Comm(M)$ when $M = \rL^\infty(X) \rtimes G$ and $G \actson (X,\mu)$ satisfies Condition \ref{cond}.

\subsection*{Definition of the C$^*$-bicategory $\Hecke(G < \cG)$}

Suppose that $G < \cG$ is a Hecke pair, meaning that $g G g^{-1} \cap G$ has finite index in $G$ for all $g \in \cG$. We define the C$^*$-bicategory $\Hecke(G < \cG)$. In the case where $G \actson (X,\mu)$ satisfies Condition \ref{cond} and $M = \rL^\infty(X) \rtimes G$, we will take $\cG$ to be the commensurator of $G$ inside $\Aut(X,\mu)$. Theorem \ref{thm.equiv-cat} will provide an equivalence of bicategories between $\Comm(M)$ and $\Hecke(G < \cG)$.

The objects (or $0$-cells) of $\Hecke(G < \cG)$ are the pairs $(G_1,\Omega_1)$ where $G_1 < \cG$ is commensurate with $G$ and $\Omega_1 \in \rZ^2(G_1,\T)$ is the obstruction $2$-cocycle of a finite dimensional projective representation.

The \emph{morphisms (or $1$-cells)} from $(G_1,\Omega_1)$ to $(G_2,\Omega_2)$ are called correspondences and defined as the set of triplets $C = (I,\cK,\pi)$ where $I \subset \cG$ is a subset that is the union of finitely many cosets $G_1 y$ and also the union of finitely many cosets $x G_2$, and if
\begin{itemize}
\item $(\cK_x)_{x \in I}$ is a family of finite dimensional Hilbert spaces,
\item $\pi(g,x,h) : \cK_x \recht \cK_{gxh}$ are unitary operators satisfying
$$\pi(g',gxh,h') \circ \pi(g,x,h) = \Omega_1(g',g) \; \pi(gg',x,hh') \; \Omega_2(h,h') \; .$$
\end{itemize}

The composition or tensor product of correspondences is defined as follows. Assume that $C = (I,\cK,\pi)$ is a correspondence from $(G_1,\Omega_1)$ to $(G_2,\Omega_2)$ and that $C' = (I',\cK',\pi')$ is a correspondence from $(G_2,\Omega_2)$ to $(G_3,\Omega_3)$. We define the tensor product $C\dpr = C \ot C'$ as $C\dpr = (I\dpr,\cK\dpr,\pi\dpr)$ given by the following formulae.

First define the set $I \times_{G_2} I'$ as the set of orbits for the action $k \cdot (x,y) = (xk^{-1},ky)$ of $G_2$ on $I \times I'$. For every orbit $G_2 \cdot (x,y) \in I \times_{G_2} I'$, define\footnote{More canonically, $\cL_{G_2\cdot(x,y)}$ is defined as the vector space consisting of families of vectors $(\xi_z)_{z \in G_2\cdot(x,y)}$ where $\xi_{k\cdot (x,y)} \in \cK_{xk^{-1}} \ot \cK'_{ky}$ and $\xi_{k\cdot(x,y)} = (\pi(e,x,k^{-1}) \ot \pi'(k,x,e))\xi_e$ for all $k \in G_2$.} the finite dimensional Hilbert space $\cL_{G_2 \cdot (x,y)} = \cK_x \ot \cK'_y$. Put $I\dpr := II'$.
The map $\theta : I \times_{G_2} I' \recht I\dpr : G_2 \cdot (x,y) \mapsto xy$ is finite-to-one. So, for every $r \in I\dpr$, we define the finite dimensional Hilbert space $\cK\dpr_r$ as the direct sum of all $\cL_{s}$, $s \in \theta^{-1}(r)$. Under the appropriate identifications, we define $\pi\dpr(g,r,h)$ as the direct sum of the operators $\pi(g,x,e) \ot \pi(e,y,h)$ when $x \in I, y \in I', xy = r$.

We call $T$ a \emph{morphism (or $2$-cell)} between the correspondences $(I,\cK,\pi)$ and $(I',\cK',\pi')$ if for every $x \in I \cap I'$, we have $T_x \in \B(\cK'_x,\cK_x)$ such that $T_{gxh} \pi'(g,x,h) = \pi(g,x,h) T_x$ for all $g \in G_1, x \in I \cap I', h \in G_2$.

One checks that $\Hecke(G < \cG)$ is a C$^*$-bicategory.

\subsection*{Bimodule functor from $\Hecke(G < \cG)$ to $\Comm(M)$}

Assume that $G \actson (X,\mu)$ is a free, weakly mixing, p.m.p.\ action. Define $\cG$ as the commensurator of $G$ inside $\Aut(X,\mu)$. By \cite[Lemma 6.11]{Va07}, $\cG$ acts freely on $(X,\mu)$. Put $A = \rL^\infty(X,\mu)$ and $M := A \rtimes G$. We define the bifunctor $\Bim$ from $\Hecke(G < \cG)$ to $\Comm(M)$.

Whenever $C = (I,\cK,\pi)$ is a correspondence between $(G_1,\Omega_1)$ and $(G_2,\Omega_2)$, we define the\linebreak $(A \rtimes_{\Omega_1} G_1)$-$(A \rtimes_{\Omega_2} G_2)$-bimodule $\Bim(C)$ given by
\begin{equation}\label{eq.elem-bim}
\begin{split}
\cH &= \bigoplus_{x \in I} (\cK_x \ot \rL^2(X,\mu)) \quad\text{and} \\  (au_g) \cdot (\xi \ot d) \cdot (b u_h) & = \pi(g,x,h)\xi \ot a \si_g(d) \si_{gx}(b) \quad\text{for all}\;\; \xi \in \cK_x \; .
\end{split}
\end{equation}
The left and right dimensions of the bimodule $\Bim(C)$ are given by
$$\diml(\Bim(C)) = \sum_{x \in G_1 \backslash I} \dim(\cK_x) \quad\text{and}\quad \dimr(\Bim(C)) = \sum_{x \in I / G_2} \dim(\cK_x) \; .$$

Whenever $T$ is a morphism between the correspondences $(I,\cK,\pi)$ and $(I',\cK',\pi')$, we define the bimodular operator $\Bim(T) : \Bim(I',\cK',\pi') \recht \Bim(I,\cK,\pi)$ defined by $\Bim(T) := \oplus_{x \in I \cap I'} (T_x \ot 1)$. In this formula, it is understood that $\Bim(T)$ is zero on $\cK_x \ot \rL^2(X,\mu)$ when $x \in I' \setminus I$.

\begin{proposition}
Suppose that $G \actson (X,\mu)$ is a free, weakly mixing, p.m.p.\ action and that $\cG$ denotes the commensurator of $G$ inside $\Aut(X,\mu)$. Put $M = \rL^\infty(X) \rtimes G$.

Then, $\Bim$ is a bifunctor from $\Hecke(G < \cG)$ to $\Bim(M)$. Moreover, $\Bim$ is isomorphic on the level of $2$-cells. More concretely, $\Bim$ defines a bijective isomorphism between the morphisms $(I,\cK,\pi) \recht (I',\cK',\pi')$ and the bounded bimodular operators $\Bim(I,\cK,\pi) \recht \Bim(I',\cK',\pi')$.
\end{proposition}
\begin{proof}
It is straightforward to check that $\Bim$ is a bifunctor. So, assume that $S$ is a bounded bimodular operator from $\Bim(I,\cK,\pi)$ to $\Bim(I',\cK',\pi')$. In particular, for all $x \in I$ and $y \in I'$, $S_{y,x} : \cK_x \ot \rL^2(X,\mu) \recht \cK'_y \ot \rL^2(X,\mu)$ is a bounded operator satisfying
$$S_{y,x}\bigl((1 \ot a) \xi (1 \ot \si_x(b))\bigr) = (1 \ot a) \, S_{y,x}(\xi) \, (1 \ot \si_y(b))$$
for all $a,b \in A := \rL^\infty(X)$. The commutation with $a \in A$ implies that $S_{y,x} \in \B(\cK_x,\cK_y) \ot A$. Since $\cG$ acts freely on $(X,\mu)$, the commutation with $b \in A$ then forces $x=y$, unless $S_{y,x} = 0$. We conclude that $S$ is the direct sum of operators $S_x \in \B(\cK_x) \ot A$, $x \in I \cap I'$.

Observe that the finite dimensional subspace $\cK_x \ot 1 \subset \cK_x \ot \rL^2(X)$ is globally invariant under the unitary operators $\xi \mapsto u_g \xi u_{x^{-1}g x}^*$, $g \in G_1 \cap x G_2 x^{-1}$. By weak mixing, it follows that $S_x(\cK_x \ot 1) \subset \cK_x \ot 1$. So, $S_x = T_x \ot 1$ for some $T_x \in \B(\cK_x)$.

We have found a morphism $T$ from $(I,\cK,\pi)$ to $(I',\cK',\pi')$ such that $S = \Bim(T)$.
\end{proof}

\subsection*{Equivalence of the bicategories $\Hecke(G < \cG)$ and $\Comm(M)$}

We then arrive at a categorical reformulation of Theorem \ref{thm.main}.

\begin{theorem}\label{thm.equiv-cat}
Let $G \actson (X,\mu)$ satisfy Condition \ref{cond}. Put $M = \rL^\infty(X) \rtimes G$ and denote by $\cG$ the commensurator of $G$ inside $\Aut(X,\mu)$. Then, $\Bim$ is an equivalence between the bicategories $\Hecke(G < \cG)$ and $\Comm(G)$.

More concretely,
\begin{itemize}
\item the II$_1$ factors that are commensurable with $M$ are, up to stable isomorphism, precisely the II$_1$ factors $A \rtimes_{\Omega_1} G_1$, where $(G_1,\Om_1)$ is a $0$-cell of $\Hecke(G < \cG)$;
\item the finite index $(A \rtimes_{\Omega_1} G_1)-(A \rtimes_{\Omega_2} G_2)$-bimodules are, up to unitary equivalence, precisely the bimodules $\Bim(C)$, where $C$ is a $1$-cell of $\Hecke(G < \cG)$;
\item the bounded bimodular operators between $\Bim(C)$ and $\Bim(C')$ are precisely the operators $\Bim(T)$, where $T$ is a $2$-cell of $\Hecke(G < \cG)$.
\end{itemize}
\end{theorem}

\begin{proof}
Let $C=(I,\cK,\pi)$ be a correspondence between $(G_1,\Omega_1)$ and $(G_2,\Omega_2)$.
It suffices to clarify the relation between the bimodules $\Bim(C)$ and the bimodules $\cK(\gamma,\pi)$ appearing in the formulation of Theorem \ref{thm.main}.

Put $P = \rL^\infty(X) \rtimes_{\Om_1} G_1$ and $Q = \rL^\infty(X) \rtimes_{\Om_2} G_2$. Choose $\gamma_1,\ldots,\gamma_n \in I$ such that $I$ is the disjoint union of the $G_1 \gamma_k G_2$. Define $\cK_k := \cK_{\gamma_k}$ and define the map
$$\pi_k : G_1 \cap \gamma_k G_2 \gamma_k^{-1} \recht \cU(\cK_k) : \pi_k(g) := \overline{\Omega_2 \circ \Ad \gamma_k^{-1}}(g^{-1},g) \; \pi(g,\gamma_k,\gamma_k^{-1}g^{-1} \gamma_k) \; .$$
Then, $\pi_k$ is a projective representation with obstruction $2$-cocycle $\Om_1 \overline{\Om_2 \circ \Ad \gamma_k^{-1}}$. One checks easily that
$$\bim{P}{\Bim(I,\cK,\pi)}{Q} \cong \bigoplus_{k=1}^n \bim{P}{\cK(\gamma_k,\pi_k)}{Q} \; .$$
\end{proof}

\end{document}